\documentclass[a4paper,english,oneside]{article}
\usepackage[T1]{fontenc}
\usepackage[utf8]{inputenc} 

\usepackage{geometry}
\usepackage{amsmath}
\usepackage{amsfonts}
\usepackage{amssymb}
\usepackage{amsthm}
\usepackage[english]{babel}
\usepackage{color}
\usepackage{lmodern}
\usepackage{pstricks-add}
\usepackage{hyperref}
\hypersetup{hidelinks}

\usepackage{mathscinet}
\usepackage{enumitem}
\usepackage[cal=boondoxo,bb=ams]{mathalfa}
\usepackage{microtype}

\usepackage{mathtools}
\usepackage{esint}

\theoremstyle{plain}
\newtheorem{theorem}{Theorem}[section]
\newtheorem{lemma}[theorem]{Lemma}

\newtheorem{proposition}[theorem]{Proposition}

\theoremstyle{definition}
\newtheorem{definition}[theorem]{Definition}

\theoremstyle{remark}
\newtheorem{remark}{Remark}

    \DeclareMathOperator\tr{Tr}
     
\begin{document}

\title{Semilinear wave equation on compact Lie groups}

\author{ Alessandro Palmieri} 
\date{}


\maketitle

\begin{abstract}

In this note, we study the semilinear wave equation with power nonlinearity $|u|^p$ on compact Lie groups. First, we prove a local in time existence result in the energy space via Fourier analysis on compact Lie groups. Then, we prove a blow-up result for the semilinear Cauchy problem for any $p>1$, under suitable sign assumptions for the initial data. Furthermore, sharp lifespan estimates for local (in time) solutions are derived.
\end{abstract}

\begin{flushleft}
\end{flushleft}

\begin{flushleft}
\textbf{AMS Classification (2020)} Primary:  35B44, 35L71; Secondary: 43A30, 43A77, 58J45
\end{flushleft}

\section{Introduction}

Let $\mathbb{G}$ be a compact Lie group and let $\mathcal{L}$ be the Laplace -- Beltrami operator on $\mathbb{G}$ (which coincides with the Casimir element of the enveloping algebra). 
In the present work, we prove a blow-up result for the Cauchy problem for the semilinear wave equation with power nonlinearity, namely,
\begin{align}\label{semilinear wave CP}
\begin{cases} \partial_t^2 u-\mathcal{L} u=| u|^p, &  x\in \mathbb{G}, \ t>0,\\
u(0,x)=\varepsilon u_0(x), & x\in \mathbb{G}, \\ \partial_t u(0,x)=\varepsilon u_1(x), & x\in \mathbb{G},
\end{cases}
\end{align} where $p>1$ and $\varepsilon$ is a positive constant describing the smallness of Cauchy data. 

Throughout the paper $L^q(\mathbb{G})$ denotes the space of $q$ -- summable functions on $\mathbb{G}$ with respect to the normalized Haar measure for $1\leqslant q < \infty$ (respectively, essentially bounded for $q=\infty$) and for $s>0$ and $q\in(1,\infty)$ the Sobolev space $H^{s,q}_\mathcal{L}(\mathbb{G)}$ is defined as the space $$H^{s,q}_\mathcal{L}(\mathbb{G)}\doteq \left\{ f\in L^q(\mathbb{G}): (-\mathcal{L})^{s/2}f\in L^q(\mathbb{G})\right\}$$ endowed with the norm
$\|f\|_{H^{s,q}_\mathcal{L}(\mathbb{G)}}\doteq\| f\|_{L^q(\mathbb{G})}+\| (-\mathcal{L})^{s/2} f\|_{L^q(\mathbb{G})}.$ As customary, the Hilbert space $H^{s,2}_\mathcal{L}(\mathbb{G)}$ will be simply denoted by $H^{s}_\mathcal{L}(\mathbb{G)}$.

For the classical semilinear wave equation in $\mathbb{R}^n$ it has been proved that the critical exponent is the positive root $p_{\mathrm{Str}}(n)$ of the quadratic equation $$(n-1)p^2-(n+1)p-2=0,$$ which is the so -- called Strauss exponent, named after the author of \cite{Str81}. For a detailed overview and a complete list of references on the proof of the Strauss' conjecture we refer to the introduction of \cite{TakWak11}.

In the not -- flat case, the semilinear wave equation has been studied in many Lorentzian metrics such as  Schwarzschild spacetime \cite{CG06,LLM20}, de Sitter spacetime \cite{Yag09,YagGal09,Yag12,GalYag2017,ER18} and Einstein -- de Sitter spacetime \cite{GalYag17EdS,TW20,PTY20,P20EdS}. In the present work, we will consider the case of semilinear wave equations in compact Lie groups, which seems to be an almost unexplored topic in the literature, up to the knowledge of the author. 

By employing Fourier analysis for compact Lie groups, the local well -- posedness of the semilinear Cauchy problem \eqref{semilinear wave CP} in the energy evolution space, that is $\mathcal{C}([0,T],H^1_{\mathcal{L}}(\mathbb{G}))\cap \mathcal{C}^1([0,T],L^2(\mathbb{G}))$, will be proved. In particular, a Gagliardo -- Nirenberg type inequality recently proved in \cite{RY19} will be used in order to estimate the power nonlinearity in $L^2(\mathbb{G})$. As byproduct of this local existence result we will obtain lower bound estimates for the lifespan of a local in time solution (that is, the maximal existence time of a solution).

Finally, we will prove the nonexistence of globally in time defined  solutions to \eqref{semilinear wave CP} for any $p>1$ and regardless of the size of initial data, provided that the data fulfill suitable sign assumptions.

\subsection{Main results}

Let us get started  with the statement of the local existence result for the semilinear Cauchy problem \eqref{semilinear wave CP}.

\begin{theorem}  \label{Thm loc esistence}
Let $\mathbb{G}$ be a compact, connected Lie group. Let us assume that the topological dimension $n$ of $\mathbb{G}$ satisfies $n\geqslant 3$. Let $(u_0,u_1)\in H_\mathcal{L}^1(\mathbb{G})\times L^2(\mathbb{G})$ and let $p>1$ such that $p\leqslant \frac{n}{n-2}$. Then, there exists $T=T(\varepsilon)>0$ such that the Cauchy problem \eqref{semilinear wave CP} admits a uniquely determined mild solution $$u\in \mathcal{C}\left([0,T],H^1_\mathcal{L}(\mathbb{G})\right)\cap \mathcal{C}^1\left([0,T],L^2(\mathbb{G})\right).$$
Furthermore, the lifespan $T$ satisfies the following lower bound estimates
\begin{align} \label{lower bound lifespan}
 T(\varepsilon) \geqslant\begin{cases}
c\varepsilon^{-\frac{p-1}{p+1}} & \mbox{if} \ u_1 \neq 0, \\
c\varepsilon^{-\frac{p-1}{2}} & \mbox{if} \ u_1 = 0,
\end{cases}
\end{align}
where the positive constant $c$ is independent of $\varepsilon$.

\end{theorem}

\begin{remark} \label{Remark ub p} The upper bound assumption for the exponent $p$ in Theorem \ref{Thm loc esistence} is made in order to apply a Gagliardo -- Nirenberg type inequality  proved in \cite[Remark 1.7]{RY19}. The restriction $n\geqslant 3$ is made to fulfill the assumptions for the employment of such inequality as well and it can be avoided by looking for solutions in weaker spaces than the  one in the statement of Theorem \ref{Thm loc esistence}, i.e., in $\mathcal{C}\left([0,T],H^s_\mathcal{L}(\mathbb{G})\right)\cap \mathcal{C}^1\left([0,T],L^2(\mathbb{G})\right)$ for some $s\in (0,1)$.
\end{remark}

We introduce now a suitable notion of energy solutions  for the semilinear Cauchy problem \eqref{semilinear wave CP}.

\begin{definition} \label{Definition energy sol}
Let $(u_0,u_1)\in H_\mathcal{L}^1(\mathbb{G})\times L^2(\mathbb{G})$. We say that $$u\in \mathcal{C}\left([0,T),H^1_\mathcal{L}(\mathbb{G})\right)\cap \mathcal{C}^1\left([0,T),L^2(\mathbb{G})\right)\cap L^p_{\mathrm{loc}}\left([0,T)\times\mathbb{G}\right)$$ is an \emph{energy solution} on $[0,T)$ to \eqref{semilinear wave CP} if $u$ fulfills the integral relation
\begin{align}
& \int_{\mathbb{G}} \partial_t u(t,x) \psi(t,x)\, \mathrm{d}x -\int_{\mathbb{G}}  u(t,x) \psi_s(t,x)\, \mathrm{d}x  -\varepsilon \int_{\mathbb{G}} u_1(x) \psi(0,x)\, \mathrm{d}x  \notag \\
 & \qquad +\varepsilon \int_{\mathbb{G}}  u_0(x) \psi_s(0,x)\, \mathrm{d}x + \int_0^t\int_\mathbb{G} u(s,x) \big(\psi_{ss}(s,x)-\mathcal{L}\psi(s,x)\big) \mathrm{d}x \, \mathrm{d}s  \notag \\
 & \quad = \int_0^t\int_\mathbb{G} |u(s,x)|^p \psi(s,x)\,  \mathrm{d}x \, \mathrm{d}s \label{def energ sol int relation}
\end{align} for any $\psi\in\mathcal{C}^\infty_0([0,T)\times \mathbb{G})$ and any $t\in (0,T)$.
\end{definition}

\begin{theorem}  \label{Thm blow up}
Let $\mathbb{G}$ be a compact Lie group. Let $(u_0,u_1)\in H_\mathcal{L}^1(\mathbb{G})\times L^2(\mathbb{G})$ be nonnegative and not both trivial functions. 
 Let $p>1$ and let  $$u\in \mathcal{C}\left([0,T),H^1_\mathcal{L}(\mathbb{G})\right)\cap \mathcal{C}^1\left([0,T),L^2(\mathbb{G})\right)\cap L^p_{\mathrm{loc}}\left([0,T)\times\mathbb{G}\right)$$ be an energy solution to \eqref{semilinear wave CP} according to Definition \ref{Definition energy sol} with lifespan $T=T(\varepsilon)$.  Then, there exists a positive constant $\varepsilon_0=\varepsilon_0(u_0,u_1,p)>0$ such that for any $\varepsilon\in (0,\varepsilon_0]$ the energy solution $u$ blows up in finite time. Moreover, the upper bound estimates for the lifespan 
\begin{equation} \label{upper bound T}
T(\varepsilon)\leqslant
\begin{cases}
 C \varepsilon^{-\frac{p-1}{p+1}} & \mbox{if} \ u_1\neq 0, \\
 C \varepsilon^{-\frac{p-1}{2}} & \mbox{if} \ u_1= 0, 
 \end{cases}
\end{equation} hold, where the constant $C>0$ is independent of $\varepsilon$.
\end{theorem}


\begin{remark} Combining \eqref{lower bound lifespan} and \eqref{upper bound T}, we obtain the sharp lifespan estimates
\begin{align*}
\begin{cases}
c  \varepsilon^{-\frac{p-1}{p+1}} \leqslant T(\varepsilon)\leqslant C \varepsilon^{-\frac{p-1}{p+1}} & \mbox{if} \ u_1\neq 0, \\
 c  \varepsilon^{-\frac{p-1}{2}} \leqslant T(\varepsilon)\leqslant C \varepsilon^{-\frac{p-1}{2}} & \mbox{if} \ u_1= 0, 
 \end{cases}
\end{align*} for local in time solutions to \eqref{semilinear wave CP}. Therefore, the nontriviality of $u_1$ plays a crucial role on the lifespan estimates, analogously to what happens for the Euclidean case (cf. \cite{Tak15}).
\end{remark}

\subsection*{Notations} Hereafter we use the following notations: $\mathcal{L}$ denotes the Laplace -- Beltrami operator on $\mathbb{G}$; $\tr(A)= \sum_{j=1}^d a_{jj}$ and $A^*=(\overline{a_{ji}})_{1\leqslant i,j\leqslant d}$ denote the trace  and the adjoint matrix of $A=(a_{ij})_{1\leqslant i,j\leqslant d} \in \mathbb{C}^{d\times d}$, respectively; $I_d\in\mathbb{C}^{d\times d}$ denotes the identity matrix; $\mathrm{d}x$ stands for the normalized Haar measure on the compact group $\mathbb{G}$; finally, when there exists a positive constant $C$ such that $f\leqslant Cg$ we write $f\lesssim g$. 

\section{Local  existence in energy space} \label{Section local existence}

In this section, we will prove Theorem \ref{Thm loc esistence}. First, we recall the notion of mild solutions to \eqref{semilinear wave CP}. By Duhamel's principle, the solution to the linear problem
\begin{align}\label{linear CP wave}
\begin{cases} \partial_t^2 u-\mathcal{L} u =F(t,x), &  x\in \mathbb{G}, \ t>0,\\
u(0,x)= u_0(x), & x\in \mathbb{G}, \\ \partial_t u(0,x)= u_1(x), & x\in \mathbb{G}
\end{cases}
\end{align}  can be represented as
\begin{align*}
u(t,x) = u_0(x)\ast_{(x)} E_0(t,x) + u_1(x)\ast_{(x)} E_1(t,x)+ \int_0^t F(s,x)\ast_{(x)} E_1(t-s,x) \, \mathrm{d}s,
\end{align*}
 where $E_0(t,x)$ and $E_1(t,x)$ denote, respectively, the fundamental solutions to \eqref{linear CP wave} in the homogeneous case $F= 0$ with initial data $(u_0,u_1)=(\delta_0,0)$ and $(u_0,u_1)=(0,\delta_0)$.
Let us point out that, in order to get the previous representation formula, we applied the invariance by time translations for the wave  operator $\partial_t^2-\mathcal{L}$ and the property $L\big(v\ast_{(x)}E_1(t,\cdot)\big)=v\ast_{(x)}L(E_1(t,\cdot))$ for any left -- invariant differential operator $L$ on $\mathbb{G}$.

Also, $u$ is said a \emph{mild solution} to \eqref{semilinear wave CP} on $[0,T]$ if $u$ is a fixed point for the nonlinear integral operator $N$ defined by
\begin{align*}
N: u\in X(T) \to N u(t,x) & \doteq \varepsilon\, u_0(x)\ast_{(x)} E_0(t,x) +\varepsilon \,
 u_1(x)\ast_{(x)} E_1(t,x) + \int_0^t |u(s,x)|^p\ast_{(x)} E_1(t-s,x) \, \mathrm{d}s
\end{align*} on the evolution space $X(T)\doteq\mathcal{C}\left([0,T],H^1_\mathcal{L}(\mathbb{G})\right)\cap \mathcal{C}^1\left([0,T],L^2(\mathbb{G})\right)$, endowed with the norm
$$\|u\|_{X(T)}\doteq \sup_{t\in[0,T]} \left(\alpha(t)^{-1}\|u(t,\cdot)\|_{L^2(\mathbb{G})}+\|(-\mathcal{L})^{1/2}u(t,\cdot)\|_{L^2(\mathbb{G})}+\|\partial_t u(t,\cdot)\|_{L^2(\mathbb{G})}\right),$$ where 
\begin{align*}
\alpha (t) \doteq 
\begin{cases} (1+t) &  \mbox{if} \ u_1\neq 0; 
\\ 1 & \mbox{if} \ u_1=0.
\end{cases}
\end{align*} As we will see in Proposition \ref{Prop L^2-L^2 estimates}, the function $\alpha$ represents the long time behavior of the $L^2(\mathbb{G})$ norm of the solution to the corresponding linear homogeneous Cauchy problem. The introduction of this time -- dependent weight in the norm on $X(T)$ will allow us to determine sharp lower bound estimates for the lifespan of local in time solutions to \eqref{semilinear wave CP}.

By employing Banach's fixed point theorem, we will show that $N$ admits a uniquely determined fixed point and we will provide a lower bound estimate for the lifespan $T=T(\varepsilon)>0$. However, before considering the semilinear Cauchy problem, we determine $L^2(\mathbb{G})$ -- $L^2(\mathbb{G})$ estimates for the solution to the corresponding linear homogeneous problem via the group Fourier transform with respect to the spatial variable. After that these estimates will have been obtained, we could show the local well -- posedness for \eqref{semilinear wave CP}  by applying a Gagliardo -- Nirenberg type inequality derived recently in \cite{RY19} (cf. Lemma \ref{Lemma GN ineq L2}). The remaining part of this section is organized in the following way: in Section \ref{Subsection GFT} we recall the main tools from Fourier Analysis on compact Lie groups which are necessary for our approach; hence, in Section \ref{Subsection L^2-L^2 est} $L^2(\mathbb{G})$ -- $L^2(\mathbb{G})$ estimates for the solution of the corresponding homogeneous linear problem and its first order derivatives are derived; finally, in Section \ref{Subsection Fixed Point} it will be shown that the operator $N$ is a contraction on $X(T)$.

\subsection{Group Fourier transform} \label{Subsection GFT}

Let us recall some results on Fourier Analysis on compact Lie groups as in \cite[Section 2.1]{Pal20WEd}. For further details on this topic we refer to the monographs \cite{RT10,FR16}. 

A \emph{continuous unitary representation} $\xi:\mathbb{G}\to\mathbb{C}^{d_\xi\times d_\xi}$ of dimension $d_\xi$ is a continuous group homomorphism from $\mathbb{G}$ to the group of unitary matrix $\mathrm{U}(d_\xi,\mathbb{C})$, that is, $\xi(xy)=\xi(x)\xi(y)$ and $\xi(x)^*=\xi(x)^{-1}$  for all $x,y\in\mathbb{G}$ and the elements $\xi_{ij}:\mathbb{G}\to \mathbb{C}$ of the matrix representation $\xi$ are continuous functions for all $i,j\in \{1,\ldots,d_\xi\}$. Two representations $\xi,\eta$ of $\mathbb{G}$ are said \emph{equivalent} if there exists an invertible intertwining operator $A$ such that $A \xi(x) = \eta(x) A$ for any $x\in \mathbb{G}$. A subspace $W\subset \mathbb{C}^{d_\xi}$ is $\xi$ -- invariant if $\xi(x) \cdot W\subset W$ for any $x\in \mathbb{G}$. A representation $\xi$ is \emph{irreducible} if the only $\xi$ -- invariant subspaces are the trivial ones $\{0\},\mathbb{C}^{d_\xi}$.

The unitary dual $\widehat{\mathbb{G}}$ of the compact Lie group $\mathbb{G}$ consists of the equivalence class $[\xi]$ of continuous irreducible unitary representation $\xi:\mathbb{G}\to\mathbb{C}^{d_\xi\times d_\xi}$. 

Given $f\in L^1(\mathbb{G})$, its Fourier coefficients at $[\xi]\in\widehat{\mathbb{G}}$ is defined by
\begin{align*}
\widehat{f}(\xi)\doteq \int_{\mathbb{G}} f(x) \xi(x)^*  \mathrm{d}x \in \mathbb{C}^{d_\xi\times d_\xi},
\end{align*} where the integral is taken with respect to the Haar measure on $\mathbb{G}$.

If $f\in L^2(\mathbb{G})$, then, the Fourier series representation for $f$ is given by
$$f(x) = \sum_{[\xi]\in\widehat{\mathbb{G}}} d_\xi \tr\big(\xi(x)\widehat{f}\ \big),$$
 where hereafter just one irreducible unitary matrix representation is picked in the sum for each equivalence class $[\xi]$ in $\widehat{\mathbb{G}}$.
Furthermore, for $f\in L^2(\mathbb{G})$ Plancherel formula takes the following form
\begin{align} \label{Plancherel formula G}
\| f\|^2_{L^2(\mathbb{G})} =  \sum_{[\xi]\in\widehat{\mathbb{G}}} d_\xi \big\| \widehat{f}(\xi)\big\|^2_{\mathrm{HS}}, 
\end{align}  where the Hilbert -- Schmidt norm of the matrix $\widehat{f}(\xi)$ is defined as follows:
\begin{align*}
 \big\| \widehat{f}(\xi)\big\|^2_{\mathrm{HS}}  \doteq  \tr\big( \widehat{f}(\xi) \widehat{f}(\xi)^*\big) = \sum_{i,j=1}^{d_\xi} \big| \widehat{f}(\xi)_{ij}\big|^2.
\end{align*}

For our analysis it is very important to understand the behavior of the group Fourier transform  with respect to the Laplace -- Beltrami operator $\mathcal{L}$. Given $[\xi]\in\widehat{\mathbb{G}}$, then, all $\xi_{ij}$ are eigenfunctions for $\mathcal{L}$ with the same not positive eigenvalue $-\lambda_\xi^2$, namely,
\begin{align*}
-\mathcal{L} \xi_{ij}(x) = \lambda^2_\xi \, \xi_{ij}(x) \quad\mbox{for any} \ x\in\mathbb{G} \ \mbox{and for all} \ i,j\in \{1,\ldots,d_\xi\}.
\end{align*} This means that the symbol of $\mathcal{L}$ is 
\begin{align} \label{symbol Laplace-Beltrami}
\sigma_\mathcal{L} (\xi)= - \lambda^2_\xi I_{d_\xi},
\end{align} that is, $\widehat{\mathcal{L} f}(\xi)= \sigma_\mathcal{L} (\xi) \widehat{f}(\xi) =- \lambda^2_\xi \widehat{f}(\xi)$ for any $[\xi]\in\widehat{\mathbb{G}}$.

Finally, through Plancherel formula for $s>0$ we have
\begin{align*}
\|f\|^2_{\dot{H}^s_\mathcal{L}(\mathbb{G})} =\| (-\mathcal{L})^{s/2}f\|^2_{L^2(\mathbb{G})} =  \sum_{[\xi]\in\widehat{\mathbb{G}}} d_\xi \lambda_\xi^{2s}\big\| \widehat{f}(\xi)\big\|^2_{\mathrm{HS}}.
\end{align*}

\subsection{$L^2(\mathbb{G})$ -- $L^2(\mathbb{G})$ estimates for the solution to the linear homogeneous problem} \label{Subsection L^2-L^2 est}

In this section, we derive $L^2(\mathbb{G})$ -- $L^2(\mathbb{G})$ estimates for the solution to the homogeneous problem associated to \eqref{semilinear wave CP}, that is,
\begin{align} \label{linear CP wave hom}
\begin{cases} \partial_t^2 u-\mathcal{L} u =0, &  x\in \mathbb{G}, \ t>0,\\
u(0,x)= u_0(x), & x\in \mathbb{G}, \\ \partial_t u(0,x)= u_1(x), & x\in \mathbb{G}.
\end{cases}
\end{align}
 We follow the main ideas from \cite{GR15,Pal20WEd}, so, the group Fourier transform with respect to the spatial variable $x$ is applied together with Plancherel identity to determine an explicit expression for the $L^2(\mathbb{G})$ norms of $u(t,\cdot)$, $(-\mathcal{L})^{1/2}u(t,\cdot)$ and $\partial_t u(t,\cdot)$, respectively.

Let $u$ be a solution to \eqref{linear CP wave hom}. Let $\widehat{u}(t,\xi) = (\widehat{u}(t,\xi)_{k\ell})_{1\leqslant k,\ell\leqslant d_\xi}\in \mathbb{C}^{d_\xi\times d_\xi}$, $[\xi]\in\widehat{\mathbb{G}}$ denote the group Fourier transform of $u$ with respect to the $x$ -- variable. Also, $\widehat{u}(t,\xi)$ is a solution of the Cauchy problem for the system of ODEs (with size depending on the representation $\xi$)
\begin{align*}
\begin{cases}
\partial_t^2 \widehat{u}(t,\xi) -\sigma_\mathcal{L}(\xi) \widehat{u}(t,\xi) =0, & t>0, \\
 \widehat{u}(0,\xi) = \widehat{u}_0(\xi), \\
 \partial_t \widehat{u}(0,\xi) =  \widehat{u}_1(\xi).
\end{cases}
\end{align*} By \eqref{symbol Laplace-Beltrami}, it follows that the previous system is decoupled in $d_\xi^2$ independent scalar ODEs, namely,
\begin{align}\label{scalar dec ODE}
\begin{cases}
\partial_t^2 \widehat{u}(t,\xi)_{k\ell} + \lambda_\xi^2 \widehat{u}(t,\xi)_{k\ell} =0, & t>0, \\
 \widehat{u}(0,\xi)_{k\ell} = \widehat{u}_0(\xi)_{k\ell}, \\
 \partial_t \widehat{u}(0,\xi)_{k\ell} =  \widehat{u}_1(\xi)_{k\ell},
\end{cases}
\end{align} for any $k,\ell\in \{1,\ldots,d_\xi\}$.
Straightforward computations provide the following representation formula for the solution to the parameter dependent ordinary linear homogeneous Cauchy problem \eqref{scalar dec ODE}
\begin{align}\label{representation u hat kl}
  \widehat{u}(t,\xi)_{k\ell} =  G_0(t,\xi)  \widehat{u}_0(\xi)_{k\ell} +G_1(t,\xi) \widehat{u}_1(\xi)_{k\ell},
\end{align} where
\begin{equation} \label{def G0 G1}
\begin{split}
 G_0(t,\xi)  & \doteq \begin{cases} \cos (\lambda_\xi t)  & \mbox{if} \ \lambda_\xi^2>0 , \\
 1  & \mbox{if} \ \lambda_\xi^2=0,   \end{cases}  \qquad
 G_1(t,\xi)  \doteq \begin{cases} \dfrac{\sin(\lambda_\xi t)}{\lambda_\xi }  & \mbox{if} \ \lambda_\xi^2>0 , \\
 t  & \mbox{if} \ \lambda_\xi^2=0. \end{cases}
 \end{split}
\end{equation} Notice that $G_0(t,\xi)=\partial_t G_1(t,\xi)$ for any $[\xi]\in\widehat{\mathbb{G}}$. We underline explicitly that $0$ is an eigenvalue for the continuous irreducible unitary representation $1:x\in\mathbb{G}\to 1\in \mathbb{C}$.

\subsubsection*{Estimate for $\|u(t)\|_{L^2(\mathbb{G})}$}
 
By \eqref{representation u hat kl} it follows 
\begin{align} \label{estimate u hat kl}|\widehat{u}(t,\xi)_{k\ell}| \leqslant |\widehat{u}_0(\xi)_{k\ell}| + t \, |\widehat{u}_1(\xi)_{k\ell}| \qquad \mbox{for any} \ t\geqslant 0.
\end{align}
Consequently, by using Plancherel formula twice, we get
\begin{align}
\| u(t,\cdot)\|^2_{L^2(\mathbb{G)}} = \sum_{[\xi]\in\widehat{\mathbb{G}}} d_\xi \sum_{k,\ell=1}^{d_\xi} |\widehat{u}(t,\xi)_{k\ell}|^2 & \lesssim  \sum_{[\xi]\in\widehat{\mathbb{G}}} d_\xi \sum_{k,\ell=1}^{d_\xi} \left( |\widehat{u}_0(\xi)_{k\ell}|^2 +t^2|\widehat{u}_1(\xi)_{k\ell}|^2 \right) = \| u_0\|^2_{L^2(\mathbb{G)}} +t^2\| u_1\|^2_{L^2(\mathbb{G)}}. \label{proof L^2 est u(t)}
\end{align}

\subsubsection*{Estimate for $\|(-\mathcal{L})^{1/2}u(t)\|_{L^2(\mathbb{G})}$}

By Plancherel formula, we find
\begin{align*}
\| (-\mathcal{L})^{1/2} u(t,\cdot)\|^2_{L^2(\mathbb{G})} = \sum_{[\xi]\in\widehat{\mathbb{G}}} d_\xi \|\sigma_{(-\mathcal{L})^{1/2} }(\xi) \widehat{u}(t,\xi)\|^2_{\mathrm{HS}} = \sum_{[\xi]\in\widehat{\mathbb{G}}} d_\xi \sum_{k,\ell=1}^{d_\xi} \lambda_\xi^2 |\widehat{u}(t,\xi)_{k\ell}|^2.
\end{align*} Since 
\begin{align}\label{estimate lambda u hat kl}
\lambda_\xi |\widehat{u}(t,\xi)_{k\ell}| \leqslant \lambda_\xi  |\widehat{u}_0(\xi)_{k\ell}| +|\widehat{u}_1(\xi)_{k\ell}|,
\end{align}  then, from the previous identity we have
\begin{align*}
\| (-\mathcal{L})^{1/2} u(t,\cdot)\|^2_{L^2(\mathbb{G})} 
& \lesssim  \sum_{[\xi]\in \widehat{\mathbb{G}}} d_\xi \sum_{k,\ell=1}^{d_\xi} \left(\lambda_\xi^2|\widehat{u}_0(\xi)_{k\ell}|^2 +|\widehat{u}_1(\xi)_{k\ell}|^2\right)  =   \|u_0\|^2_{H^1_\mathcal{L}(\mathbb{G})}+\|u_1\|^2_{L^2(\mathbb{G})},
\end{align*}  where in the last step we applied again Plancherel formula.

\subsubsection*{Estimate for $\|\partial_t u(t)\|_{L^2(\mathbb{G})}$}

Elementary computations show that for any $[\xi]\in\widehat{\mathbb{G}}$ and any $k,\ell\in\{1,\ldots,d_\xi\}$ it holds
\begin{align*}
\partial_t \widehat{u}(t,\xi)_{k\ell} = - \lambda_\xi^2 \, G_1(t,\xi)  \widehat{u}_0(\xi)_{k\ell} +G_0(t,\xi)  \widehat{u}_1(\xi)_{k\ell},
\end{align*} where $G_0(t,\xi),G_1(t,\xi)$ are defined in \eqref{def G0 G1}. Therefore,
\begin{align*}
|\partial_t \widehat{u}(t,\xi)_{k\ell}| \leqslant \lambda_\xi  |\widehat{u}_0(\xi)_{k\ell}| +|\widehat{u}_1(\xi)_{k\ell}|.
\end{align*}  Since the right -- hand side of the previous inequality is the same one as in \eqref{estimate lambda u hat kl}, we get immediately
\begin{align*}
\| \partial_t u(t,\cdot)\|^2_{L^2(\mathbb{G})} 
& \lesssim    \|u_0\|^2_{H^1_\mathcal{L}(\mathbb{G})}+\|u_1\|^2_{L^2(\mathbb{G})}.
\end{align*}

Summarizing, in this section we proved the following proposition.
\begin{proposition} \label{Prop L^2-L^2 estimates}
Let us assume $(u_0,u_1)\in H^1_\mathcal{L}(\mathbb{G})\times L^2(\mathbb{G})$ and let $$u\in \mathcal{C}\big([0,\infty),H^1_\mathcal{L}(\mathbb{G})\big)\cap  \mathcal{C}^1\big([0,\infty),L^2(\mathbb{G})\big)$$ be the solution to the homogeneous Cauchy problem \eqref{linear CP wave hom}.
 
 Then, the following $L^2(\mathbb{G})$ -- $L^2(\mathbb{G})$ estimates are satisfied
\begin{align}
\|u(t,\cdot)\|_{L^2(\mathbb{G})} &\leqslant C \Big( \|u_0\|_{L^2(\mathbb{G})}+ t \,  \|u_1\|_{L^2(\mathbb{G})}\Big), \label{L^2 norm u(t)} \\
\|(-\mathcal{L})^{1/2}u(t,\cdot)\|_{L^2(\mathbb{G})} &\leqslant C \left( \|u_0\|_{H^1_\mathcal{L}(\mathbb{G})}+\|u_1\|_{L^2(\mathbb{G})}\right), \label{L^2 norm (-L)^1/2 u(t)} \\
\| \partial_t u(t,\cdot)\|_{L^2(\mathbb{G})} &\leqslant C \left( \|u_0\|_{H^1_\mathcal{L}(\mathbb{G})}+\|u_1\|_{L^2(\mathbb{G})}\right), \label{L^2 norm ut(t)}
\end{align} for any $t\geqslant 0$, where $C$ is a positive multiplicative constant.
\end{proposition}

\subsection{Proof of Theorem \ref{Thm loc esistence}} \label{Subsection Fixed Point}

A fundamental tool to prove the local existence result is the following Gagliardo -- Nirenberg type inequality, whose proof can be found in \cite{RY19} (see also \cite[Corollary 2.3]{Pal20WEd}).

\begin{lemma}\label{Lemma GN ineq L2}
Let $\mathbb{G}$ be a connected unimodular Lie group with topological dimension $n\geqslant 3$. For any $q\geqslant 2$ such that $q\leqslant \frac{2n}{n-2}$ the following Gagliardo -- Nirenberg type inequality holds
\begin{align} \label{GN ineq L2}
\| f\|_{L^q(\mathbb{G})}\lesssim \| f\|^{\theta(n,q)}_{H^{1}_\mathcal{L}(\mathbb{G})} \| f\|^{1-\theta(n,q)}_{L^{2}(\mathbb{G})}
\end{align}
for any $f\in H^{1}_\mathcal{L}(\mathbb{G})$, where $ \theta(n,q)\doteq  n\left(\frac{1}{2}-\frac{1}{q}\right)$.
\end{lemma}

We can now proceed with the proof of Theorem \ref{Thm loc esistence}.

Let us estimate $\|Nu\|_{X(T)}$ for $u\in X(T)$. We rewrite $Nu=  u^{\mathrm{ln}}+ J u$, where 
\begin{align*}
 u^{\mathrm{ln}}(t,x) &\doteq \varepsilon\, u_0(x)\ast_{(x)} E_0(t,x) +\varepsilon\,  u_1(x)\ast_{(x)} E_1(t,x),  \\  J u(t,x) & \doteq \int_0^t |u(s,x)|^p\ast_{(x)} E_1(t-s,x) \, \mathrm{d}s.
\end{align*} 
 By Proposition \ref{Prop L^2-L^2 estimates} we get immediately $\|u^{\mathrm{ln}}\|_{X(T)}\lesssim  \varepsilon\, \|(u_0,u_1)\|_{H^1_\mathcal{L}(\mathbb{G})\times L^2(\mathbb{G})}$. On the other hand, due to the invariance by time translations of \eqref{linear CP wave hom}, it results 
\begin{align}
\|\partial_t^j (-\mathcal{L})^{i/2} Ju(t,\cdot)\|_{L^2(\mathbb{G})} & \lesssim \int_0^t (t-s)^{1-(j+i)} \| u(s,\cdot)\|^p_{L^{2p}(\mathbb{G})} \, \mathrm{d}s \notag \\ 
& \lesssim \int_0^t (t-s)^{1-(j+i)} \| u(s,\cdot)\|^{p\theta(n,2p)}_{H^{1}_{\mathcal{L}}(\mathbb{G})} \| u(s,\cdot)\|^{p(1-\theta(n,2p))}_{L^2(\mathbb{G})} \, \mathrm{d}s \notag \\
& \lesssim \int_0^t (t-s)^{1-(j+i)} \alpha(s)^p\| u\|^{p}_{X(s)}  \, \mathrm{d}s \notag\\ & \lesssim t^{2-(j+i)} \alpha(t)^p \, \| u\|_{X(t)}^p \label{estimate Ju in X(T)}
\end{align} for $i,j\in\{0,1\}$ such that $0\leqslant i+j\leqslant 1$. We stress that employing \eqref{GN ineq L2}  in the above estimate, we have to require the condition $p\leqslant \frac{n}{n-2}$ on $p$ in Theorem \ref{Thm loc esistence}. Similarly, combining H\"older's inequality and \eqref{GN ineq L2}, for $i,j\in\{0,1\}$ such that $0\leqslant i+j\leqslant 1$ we find
\begin{align}
\|\partial_t^j & (-\mathcal{L})^{i/2}   (Ju(t,\cdot)-Jv(t,\cdot))\|_{L^2(\mathbb{G})}   \notag \\ &
 \lesssim \int_0^t(t-s)^{1-(j+i)} \| |u(s,\cdot)|^p-|v(s,\cdot)|^p\|_{L^{2}(\mathbb{G})} \, \mathrm{d}s \notag \\   & \lesssim \int_0^t (t-s)^{1-(j+i)} \| u(s,\cdot)-v(s,\cdot)\|_{L^{2p}(\mathbb{G})} \left(\|u(s,\cdot)\|^{p-1}_{L^{2p}(\mathbb{G})}+ \|v(s,\cdot)\|^{p-1}_{L^{2p}(\mathbb{G})}\right) \, \mathrm{d}s \notag \\ 
& \lesssim t^{2-(j+i)} \alpha(t)^p \, \| u-v\|_{X(t)}\left(\|u\|^{p-1}_{X(t)}+\|v\|^{p-1}_{X(t)}\right). \label{estimate Ju -Jv in X(T)}
\end{align} Summarizing, we just proved 
\begin{align*}
\| N u\|_{X(t)} &\leqslant \widetilde{C}  \varepsilon\, \|(u_0,u_1)\|_{H^1_\mathcal{L}(\mathbb{G})\times L^2(\mathbb{G})} +\widetilde{C} (1+t)^{\beta(p)} \|u\|^p_{X(t)}, \\
\| N u -Nv\|_{X(t)} &\leqslant \widetilde{C} (1+t)^{\beta(p)} \| u-v\|_{X(t)}\left(\|u\|^{p-1}_{X(T)}+\|v\|^{p-1}_{X(T)}\right),
\end{align*} where $$\beta(p) \doteq \begin{cases} p+1 & \mbox{if} \ u_1\neq 0,\\  2 & \mbox{if} \ u_1= 0.\end{cases} $$ 
Consequently,  denoting $R_0=\|(u_0,u_1)\|_{H^1_\mathcal{L}(\mathbb{G})\times L^2(\mathbb{G})} $, we find that for any $R\geqslant 2\widetilde{C} R_0$ and any $t\leqslant C\varepsilon^{-\frac{p-1}{\beta(p)}}$, where $C\doteq  (4\widetilde{C} R)^{-\frac{1}{\beta(p)}}$, the mapping $N$ satisfies
\begin{align*}
\| N u\|_{X(t)} \leqslant R\varepsilon, \qquad
\| N u -Nv\|_{X(t)} \leqslant \tfrac 12 \| u-v\|_{X(t)},
\end{align*} for all $u,v\in \mathfrak{B}(R\varepsilon)\doteq \{u\in X(t): \| u\|_{X(t)}\leqslant R\varepsilon\}$, that is, $N$ is a contraction mapping on the ball $ \mathfrak{B}(R\varepsilon)$ in the Banach space $X(t)$. Thus, Banach's fixed point implies the existence of a uniquely determined fixed point $u$ for $N$, which is the mild solution to \eqref{semilinear wave CP} on $[0,t]\subset [0,T(\varepsilon)]$ we were looking for. Besides, we got the lower bound estimates $T(\varepsilon)\geqslant C\varepsilon^{-\frac{p-1}{\beta(p)}}$. This completes the proof of Theorem \ref{Thm loc esistence}.


\section{Proof of Theorem \ref{Thm blow up}}

In this section we are going to prove Theorem \ref{Thm blow up} by using a comparison argument for ordinary differential inequality of the second order (\emph{Kato's lemma}).

Let $u$ be a local in time energy solution to \eqref{semilinear wave CP} according to Definition \ref{Definition energy sol} with lifespan $T$ and let us fix $t\in (0,T)$. If we choose a cut -- off function $\psi\in \mathcal{C}^\infty_0 ([0,T)\times \mathbb{G})$ such that $\psi=1$ on $[0,t]\times\mathbb{G}$ in \eqref{def energ sol int relation}, then,
\begin{align*}
& \int_{\mathbb{G}} \partial_t u(t,x) \, \mathrm{d}x  -\varepsilon \int_{\mathbb{G}} u_1(x) \, \mathrm{d}x  = \int_0^t\int_\mathbb{G} |u(s,x)|^p   \mathrm{d}x \, \mathrm{d}s. 
\end{align*} Introducing the time -- dependent functional
\begin{align*}
U_0(t)\doteq \int_{\mathbb{G}}  u(t,x) \, \mathrm{d}x,
\end{align*}  we can rewrite the above equality in the following way
\begin{align*}
U_0'(t) - U_0'(0)=\int_0^t\int_\mathbb{G} |u(s,x)|^p   \mathrm{d}x \, \mathrm{d}s.
\end{align*} Then, $U'_0$ is differentiable with respect to $t$ and 
\begin{align*}
U_0''(t) = \int_\mathbb{G} |u(t,x)|^p   \mathrm{d}x \geqslant |U_0(t)|^p,
\end{align*}
  where in the last step we applied Jensen's inequality. We remark that $$U_0(0)=\varepsilon \int_\mathbb{G} u_0(x) \, \mathrm{d}x\geqslant 0 \qquad \mbox{and} \qquad U'_0(0)=\varepsilon \int_\mathbb{G} u_1(x) \, \mathrm{d}x\geqslant 0$$ thanks to the assumptions on the initial data in the statement of Theorem \ref{Thm blow up}. Then, employing Lemmas 2.1 and 2.2 in \cite{Tak15} (improved Kato's lemma with upper bound estimate for the lifespan) to the functional $U_0$ we conclude the proof of Theorem \ref{Thm blow up}.

\section{Final remarks}
Lately, the Cauchy problem for the semilinear damped wave equation with power nonlinearity has been studied in the framework of compact Lie groups in \cite{Pal20WEd} (in this case the differential operator is the damped wave operator $\partial_t^2-\mathcal{L}+\partial_t$). In particular, in the compact case it has been proved for any exponent $p>1$ the nonexistence of global in time solution, under certain sign assumptions for the initial data. This result is consistent with the one \cite{RY18}  for the semilinar heat equation on unimodular Lie group with polynomial volume growth. Indeed, we can read this result by saying that in the compact case the critical exponent for the semilinear damped wave equation is the Fujita exponent in the 0 -- dimensional case. Therefore, rather than the topological dimension of the group $\mathbb{G}$, it is the global dimension of $\mathbb{G}$ (which is 0 for a compact Lie group) that provides the critical exponent. We refer to \cite[Section II.4]{DER03} for an overview on the growth properties of a Lie group.

A similar analysis can be done for the semilinear wave equation in \eqref{semilinear wave CP}. In the flat case, the subcritical condition $1<p<p_{\mathrm{Str}}(n)$ is equivalent to require that the quantity $\gamma(n,p)\doteq -(n-1)p^2+(n+1)p+2$ is strictly positive. In the 0 -- dimensional case, we have $\gamma(0,p)>0$ for any $p>1$. Therefore, also for the semilinear Cauchy problem associated to wave operator $\partial_t^2-\mathcal{L}$, it is possible to say that the global dimension of $\mathbb{G}$ determines the range for $p$ in the blow -- up result (in the sense that we have just explained).

Finally, we point out that not necessarily for all semilinear hyperbolic model on a compact Lie group a blow -- up result can be proved for any $p>1$. In the forthcoming paper \cite{Pal20WEdm}, it will be shown that the combined presence of a damping term and of a mass term (0th order term with a positive multiplicative constant) modifies completely the situation since the global in time existence of small data solution can be proved in the evolution energy space without requiring any additional lower bound for $p>1$.


%
%



\end{document}